\documentclass[12pt]{amsart}
\usepackage{epic,eepic}
\newlength{\baseunit}               
\newcommand{\bpf}{\noindent {\em Proof.  }}
\newcommand{\epf}{\qed \vspace{+10pt}}
\newtheorem{tm}{Theorem}
\newtheorem{pr}[tm]{Proposition}
\newtheorem{lm}[tm]{Lemma}
\newtheorem{co}[tm]{Corollary}

\newcommand{\Q}{\mathbb{Q}}
\newcommand{\com}{\mathbb{C}}

\newcommand{\De}{\Delta}

\newcommand{\proj}{\mathbb P}
\newcommand{\eff}{\mathbb F}

\newcommand{\oh}{{\mathcal{O}}}
\newcommand{\cm}{{\mathcal{M}}}
\newcommand{\cu}{{\mathcal{U}}}

\newcommand{\tC}{\tilde{C}}
\newcommand{\hC}{\hat{C}}
\newcommand{\tg}{\tilde{g}}
\newcommand{\hg}{\hat{g}}

\newcommand{\om}{\omega}

\newcommand{\de}{\delta}

\newcommand{\la}{\lambda}
\newcommand{\Pic}{\operatorname{Pic}}

\newcommand{\mbar}{\overline{M}}
\newcommand{\cmbar}{\overline{\cm}}
\newcommand{\mgs}{\cmbar_g(\proj^2,d)^*}
\newcommand{\mos}{\cmbar_1(\proj^2,d)^*}
 
\begin{document}
\pagestyle{plain}
\title{Enumerative geometry of plane curves of low genus}
\author{Ravi Vakil}
\date{Mar. 2, 1998.}
\begin{abstract}
We collect various known results (about plane curves and the moduli
space of stable maps) to derive new recursive formulas enumerating low
genus plane curves of any degree with various behaviors.  Recursive
formulas are given for the characteristic numbers of rational plane
curves, elliptic plane curves, and elliptic plane curves with fixed
complex structure.  Recursions are also given for the number of
elliptic (and rational) plane curves with various ``codimension 1''
behavior (cuspidal, tacnodal, triple pointed, etc., as well as the
geometric and arithmetic sectional genus of the Severi variety).  We
compute the latter numbers for genus 2 and 3 plane curves as well.  We rely
on results of Caporaso, Diaz, Getzler, Harris, Ran, and especially
Pandharipande.
\end{abstract}
\maketitle
\section{Introduction}
Let $\mgs$ be the component of the stack $\cmbar_g(\proj^2,d)$
generically parametrizing maps from irreducible curves.  (All stacks
will be assumed to be Deligne-Mumford stacks.)  On the universal curve
$\cu$ over $\mgs$ (with structure map $\pi$) there are two natural
divisors, the pullback $D$ of $\oh_{\proj^2}(1)$, and the relative
dualizing sheaf $\om$.  Following the notation of \cite{dh1}, let
$A=\pi_* (D^2)$, $B=
\pi_* (D \cdot \om)$, $C=\pi_* (\om^2)$, and $TL=A+B$.  Let $\De_0$ be
the divisor generically parametrizing maps from irreducible nodal
curves, and let $\De_{i,j}$ ($0<j<d$) be the divisor generically
parametrizing maps from a reducible curve, one component of genus $i$
and mapping with degree $j$, and the other of genus $g-i$ mapping with
degree $d-j$.  Let $\De = \De_0 + \sum \De_{i,j}$.  (\cite{dh1} deals with
Severi varieties, but all arguments carry over to this situation.)
Then $TL$ is the divisor class corresponding to curves tangent to a fixed
line.  Call irreducible divisors on $\mgs$ whose general map
contracts a curve {\em enumeratively meaningless}; call other divisors
{\em enumeratively meaningful}.  Call enumeratively meaningful
irreducible divisors whose general source curve is singular {\em boundary
divisors}; these are the components of $\De$.  

When $g \leq 2$, $C$ can be expressed as a sum of boundary divisors.
When $g=0$ (resp. 1), $TL - ( \frac {d-1} d )A$ (resp. $TL-A=B$) can
be expressed as a sum of boundary divisors.  By restricting this
identity to the one-dimensional family in $\mgs$ generically
corresponding to curves through $a$ general points and tangent to
$3d+g-2$ general lines, we find recursions for
characteristic numbers (when $g \leq 1$).  Recursions for the genus 0
characteristic numbers are well-known (\cite{iqd}, \cite{ek1},
\cite{ek2}).  Algorithms to determine genus 1 characteristic numbers
are known (\cite{gp2} via descendents and topological recursions;
\cite{char} by degenerations), but the formulas given here seem less
unwieldy and more suitable for computation.

In \cite{dh1} and \cite{dh2}, many divisors on the Severi variety are
expressed as linear combinations of $A$, $B$, $C$ and boundary
divisors.  (Diaz and Harris conjecture that up to torsion, any divisor
can be so expressed.)  Modulo enumeratively meaningless divisors,
their arguments carry over to $\mgs$.  Now restrict these divisors to
the one-parameter family corresponding to curves through $3d+g-2$
general points.  If $X$ is a divisor on a curve, denote its degree by
$|X|$.  When $g \leq 1$, there are simple recursions for $|A|$, $|B|$,
$|C|$, and any boundary divisor, so we get similar recursions for
enumerative ``divisor-related behavior'' (e.g. the geometric and
arithmetic sectional genus of the Severi variety, or the number of
cuspidal or triple-pointed curves, or the number of curves through
$3d-1$ general points and with 3 collinear nodes).  Some of these
formulas were known earlier.  When $g=2$ or 3, $|A|$ and $|B|$ can be
found using \cite{rinv} or \cite{ch}, and $|C|$ is simple to compute
using \cite{m} (and, if $g=3$, Graber's algorithm \cite{g} for
counting hyperelliptic plane curves).  (When $g=2$, the number $|A|$,
and possibly $|B|$, can be quickly computed by the recursions of
Belorousski and Pandharipande (\cite{bp}).  Hence these ``codimension
1'' numbers (e.g. counting cuspidal, tacnodal, or triple-pointed genus
2 or 3 curves, or computing the geometric and arithmetic sectional
genus of the Severi variety $V^{d,2}$ or $V^{d,3}$) can be computed.

The author is grateful to Tom Graber, Pasha Belorousski, Ziv Ran, and
Rahul Pandharipande for useful discussions and advice, and to Joe
Harris for first introducing him to these questions.  He also thanks
A.  Postnikov for discussing the combinatorial background to Section
\ref{covers}.  A maple program implementing all algorithms described
here is available upon request.

\subsection{Historical notes}
For a more complete historical background and references, see the
introduction to \cite{asmooth}. 

Characteristic number problems date from the last century, and were
studied extensively by Schubert, Chasles, Halphen, Zeuthen, and
others.  Much of their work is collected in \cite{s}.  

A modern study of the enumerative geometry of cubics was undertaken
successfully in the 1980's.  Among the highlights: Sacchiero and
Kleiman-Speiser independently verified Maillard and Zeuthen's results
for cuspidal and nodal plane cubics, and Kleiman and Speiser
calculated the characteristic numbers of smooth plane cubics
(\cite{ksmooth}).  Sterz (\cite{st}) and Aluffi (\cite{asmooth})
independently constructed a smooth variety of ``smooth cubics'', and
Aluffi used this variety to compute the characteristic numbers of
smooth plane cubics and other enumerative information.

The advent of the moduli space of stable maps has had tremendous
applications in enumerative algebraic geometry; as an example
pertaining to this article, Pandharipande calculated the
characteristic numbers of rational curves in $\proj^n$ in \cite{iqd},
and Ernstr\"{o}m and Kennedy showed that the characteristic numbers of
rational curves in $\proj^2$ were encoded in a ``contact cohomology
ring'' of $\proj^2$ that is the deformation of the quantum cohomology
ring (\cite{ek1}, \cite{ek2}).

\subsection{Gromov-Witten theory}
Although it isn't evident in the presentation, the main idea came from
an attempt to understand geometrically why Gromov-Witten invariants
determine gravitational descendents in genus 1 (see \cite{kk}).  This
fact should really be seen as related to a more elementary fact of
Kodaira's, that the relative dualizing sheaf of a family of elliptic
curves can be expressed as a sum of boundary divisors
(cf. (\ref{mrshowell})).  Kodaira's relation can also be used
enumeratively, by restricting it to one-parameter families, as in this
article.  For the same reason, Belorousski and Pandharipande's new
relation in $\cmbar_{2,3}$ together with Getzler's genus 2 descendent
relations (\cite{ge}) yields recursions for all $g=2$ descendent
integrals on $\proj^2$.  However, full reconstruction in $g=2$ has not
yet been shown for arbitrary spaces -- additional universal relations
are needed.  These results may be interpreted to suggest the existence
of recursive formulas for characteristic numbers of genus 2 curves,
although the recursions are likely quite messy.  (This is quite
speculative, of course.)

\section{Characteristic numbers}
We work over an algebraically closed field of characteristic 0.
Let $R_d(a,b)$ be the number of irreducible degree $d$ rational curves
through $a$ fixed general points and tangent to $b$ fixed general
lines if $a+b=3d-1$, and 0 otherwise.  Let $R_d := R_d(3d-1,0)$ be the
number with no tangency conditions.  Let $NL_d(a,b)$ be the number of
irreducible degree $d$ rational curves through $a$ fixed general
points and tangent to $b$ fixed general lines and with a node on a
fixed line if $a+b=3d-2$, and 0 otherwise.  By \cite{dh1} (1.4) and (1.5),
\begin{equation}
\label{nl} 
NL_d(a,b) = (d-1) R_d(a+1,b) - R_d(a,b+1)/2. 
\end{equation}
Let $NP(a,b)$
be the number of irreducible degree $d$ rational curves through $a$
fixed general points and tangent to $b$ fixed general lines and with a
node at a fixed point if $a+b=3d-3$, and 0 otherwise.  Let $NP_d :=
NP_d(3d-3,0)$ be the number with no tangency conditions.  Let
$E_d(a,b)$ be the number of irreducible degree $d$ elliptic curves
through $a$ fixed general points and tangent to $b$ fixed general
lines if $a+b=3d$, and 0 otherwise.  Let $E_d := E_d(3d,0)$ be the
number with no tangency conditions.

The algorithm involves six different recursions, three of them well-known 
and three quite simple:
\begin{enumerate}
\item Calculating $R_d$ using Kontsevich's recursion (\ref{Kont}).
\item Calculating $NP_d$, in essence by using Kontsevich's recursion on 
the convex rational ruled surface $\eff_1$.
\item Calculating $E_d$ using the recursion of Eguchi, Hori, and Xiong.
\item Calculating the characteristic numbers $R_d(a,b)$ using the 
characteristic numbers of lower degree curves, or curves of the same degree
with fewer tangency conditions.
\item The same thing for $NP_d(a,b)$.
\item The same thing for $E_d(a,b)$.
\end{enumerate}

\subsection{Bertini-type preliminaries}

Assume that $W$ is a variety defined over an algebraically closed
field of characteristic 0.  Consider a family of maps from nodal
curves to $\proj^2$: $$
\begin{array}{rcccl}
\cu & & \substack{{\rho}\\{ \rightarrow}} & & \proj^2 \\
& \pi \searrow & & \swarrow & \\
& & W & & 
\end{array}
$$ We say that a map has a tangent line $l \subset \proj^2$ if the
pullback of $l$ to $\cu$ contains a point with multiplicity at least 2; similar
definitions apply for flex lines and bitangents.

Let $\om$ be the relative dualizing sheaf of $\pi$, and $D = \rho^*
\oh_{\proj^2}(1)$.  Let $A = \pi_*(D^2)$ and $B = \pi_* (D \cdot \om)$ 
for convenience.  By the Kleiman-Bertini theorem (\cite{k}) applied
to $\cu$, $D$ is base-point free, and if $V$ is any irreducible
substack of $W$, a general representative of $\pi_* D^2$
intersects $V$ properly and transversely.  (Strictly speaking,
Kleiman-Bertini should be applied to $W \times PGL_2$ with group
$PGL_2$ as follows.  There is a universal curve $(\pi,1): \cu \times
PGL_2 \rightarrow W \times PGL_2$, and the universal map to $\proj^2$
is given by $(p,g) \mapsto g \circ p$.  For the sake of brevity, we
will elide this discussion when we invoke Kleiman-Bertini in the
future.)

Next, assume that $W$ is irreducible and $\cu_w$ is smooth for general
$w \in W$.  Let $L$ be the divisor on $\cu$ that is the pullback of a
general line $l$ in $\proj^2$ (so $[L] = D$).  Then $L$ has the same
dimension as $W$, its ramification divisor is in the divisor class
$(D+\om)|_L$, and the branch divisor is in class $A+B = \pi_* (D
\cdot (D+ \om))$.  

\begin{lm} 
\label{island}
If the general curve is smooth, and the general map in the 
family factors as a simply ramified multiple cover followed
by an immersion, then:
\begin{enumerate}
\item[(a)] the branch divisor is reduced, and
\item[(b)] if $V$ is any irreducible subvariety of $W$, then (for a general 
$L$) the branch divisor intersects $V$ properly.
\end{enumerate}
\end{lm}
\bpf
For part (a), we must show that the general point of any component of
the branch divisor corresponds to a map simply tangent to the line $l$
(i.e. $l$ is not a bitangent or a flex).

The general map in the family has a finite number of bitangent and
flexes.  (The image curve has a finite number of bitangents and
flexes, as the dual of a reduced curve is a reduced curve in
characteristic 0.  The only additional bitangents and flexes must
involve the simple ramification of the map from the source to the
image.  This will yield only a finite number of each.)  By a similar
argument, any particular map has at most a one-dimensional family of
bitangent lines or flex lines; call the locus with a
positive-dimensional family of such lines $B$, a proper subvariety of
$W$.  Then (for dimensional reason), the branch divisor of the
pullback of a general line $l$ to the family meets each each
component of $B$ properly.  Hence (a) follows.

Part (b) is similar, and omitted for the sake of brevity.
\epf

Hence in a one-parameter family of maps (satisfying the conditions of
the lemma), the number of curves tangent to a general fixed line is 
$|A+B|=|D \cdot (D+\om)|$.

\begin{lm}  
\label{hut}
Let $W$ be an irreducible reduced substack of $\mgs$ whose generic
member corresponds to a map from a smooth curve.  Then the subset of
$W$ corresponding to maps through a fixed general point (resp. tangent
to a fixed general line) is of pure codimension 1, each component
generically corresponds to a map from a smooth curve, and the
corresponding Weil divisor is in class $A|_W$ (resp. $(A+B)|_W$). 
\end{lm}
\bpf
The Kleiman-Bertini argument for incidence conditions is well-known
(see \cite{fp} Section 9).  We show the result for the locus $T$ in
$W$ corresponding to maps tangent to a fixed general line.  By
purity of branch locus, $T$ is pure codimension 1 in $W$.  By Lemma
\ref{island} (a), $T$ (as a Weil divisor) is in class $(A+B)|_W$.  The
irreducible components of the (proper) substack corresponding to maps
from singular curves all meet $T$ properly by Lemma \ref{island} (b),
so the general point of each component of $T$ corresponds to a map
from a smooth curve.
\epf
\begin{co} \label{enum}
For $W$ as in Lemma \ref{hut}, such that the generic map in $W$ has
trivial automorphism group, $A^a (A+B)^{\dim W - a}$ is the solution
to the enumerative problem: how many maps in $W$ pass through $a$
general points and are tangent to $(\dim W - a)$ general lines?
\end{co}

We will need to understand the divisor $TL=A+B$ on maps from nodal
curves as well.
\begin{lm}
\label{maryann}
If $W$ is an irreducible family of maps and $\cu_w$ is a curve with one node
for a general $w \in W$, then the divisor $\pi_* (D \cdot (D+\om))$ is
the divisor corresponding to where the map from the normalization is
tangent to a fixed general line $l$, plus twice the divisor
corresponding to where the node maps to $l$.  If $V$ is any
irreducible subvariety of $W$, then this divisor mets $V$ properly
(for general $l$).
\end{lm}

\bpf 
Compare the relative dualizing sheaf of the nodal curve with the
relative dualizing sheaf of the normalization.  \epf

Next, we recall relevant facts about the moduli stack of stable maps.
The stack $\cmbar_0(\proj^2,d)$ is smooth of dimension $3d-1$.  The
stack $\mos$ is the closure (in $\cmbar_1(\proj^2,d)$) of maps that
collapse no elliptic component.  It has dimension $3d$, and it is
smooth away from the divisor where an elliptic component is collapsed
(\cite{ratell} Lemma 3.13).  In particular, if $\De$ is the union of
divisors corresponding to maps from nodal curves with no collapsed
elliptic component, then $\mos$ is smooth at the generic point of each
component of $\De$.

\begin{lm} \label{ginger}
Suppose $\De$ is the locus in $\mos$ described above, or the locus in
$\cmbar_0(\proj^2,d)$ generically corresponding to maps from curves
with one node.  Fix $a$ general points and $b$ general lines, where
$a+b = \dim \De$.  Then the intersection $\De \cdot A^a TL^b$ is
equal to the number of maps where the map from the normalization
passes through the $a$ points and is tangent to the $b$ lines; plus
twice the number where the node maps to one of the $b$ lines, and the
curve passes through the $a$ points and is tangent to the remaining
$b-1$ lines; plus four times the number where the node maps to the
intersection of two of the $b$ lines, and the curve passes through the
$a$ points and is tangent to the remaining $b-2$ lines.
\end{lm}
\bpf
This follows from the fact that the condition of requiring the node to
map to a fixed general line is transverse to any subvariety (by
Kleiman-Bertini), and Lemma \ref{maryann}.
\epf

\subsection{Incidences only}
\label{theprofessor}
We begin by considering cases with no tangencies.

Clearly $R_1=1$.  There is a well-known formula (\cite{km} Claim 5.2.1 or \cite{rt}) for computing $R_d$ inductively:
\begin{equation}
\label{Kont}
R_d = \sum_{i+j=d}  i^2 j  \left( j \binom {3d-4} {3i-2} - i \binom {3d-4} {3i-1} \right) R_i R_j.
\end{equation}
One proof involves studying rational curves through $3d-2$ fixed
points, two of which are marked $p$ and $q$, and two marked points $r$
and $s$ on fixed general lines, and pulling back an equivalence on
$\Pic \cmbar_{0,4}$.  The same ``cross-ratio'' trick gives a
recursion for $NP_d$:
\begin{eqnarray}
NP_d &=& \sum_{i+j=d} (ij-1)i \left(j \binom {3d-6} {3i-3} - i \binom{3d-6}{3i-2} \right)  R_i R_j \nonumber \\
 & & +  \sum_{i+j=d} i j \left(2ij \binom {3d-6} {3i-4} - i^2 \binom{3d-6}{3i-3} - j^2 \binom {3d-6} {3i-5} \right)  NP_i R_j. \label{mrhowell}
\end{eqnarray}
(Pandharipande gives another recursion for $NP_d$ in \cite{canonical} Section 3.4.)
The Eguchi-Hori-Xiong formula (proved by Pandharipande in \cite{pgetz} using
Getzler's relation) gives $E_d$:
\begin{equation}
\label{EHX}
E_d = \frac 1 {12} \binom d 3 R_d + \sum_{i+j=d} \frac {ij (3i-2)} 9 
\binom {3d-1} {3j} R_i E_j.
\end{equation}
(Remarkably, there is still no geometric proof known of this result.)

\subsection{Swapping incidences for tangencies:  genus 0}

From \cite{iqd} Lemma 2.3.1, in $\Pic ( \cmbar_0(\proj^2,d)) \otimes \Q$,
\begin{equation}
\label{TL0}
TL = \frac {d-1} d A +  \sum_{j=0}^{[d/2]} \frac {j (d-j) } d \De_{0,j}.
\end{equation}

Apply this rational equivalence to the one-parameter family
corresponding to degree $d$ rational curves through $a$
general points and tangent to $b$ general lines (where $a+b=3d-2$) to
get:
\begin{eqnarray*}
R_d(a,b+1) &=& \frac {d-1} d R_d(a+1,b)  \\
& & + \sum_{i+j=d} \frac {ij} {2d} \Biggl[
\sum_{\substack{{a_i + a_j=a}\\{b_i+b_j=b}}} 
  \binom a {a_i} \binom b {b_i} (ij) R_i(a_i,b_i) R_j(a_j,b_j) \\
& & + 4b \sum_{\substack{{a_i + a_j=a+1}\\{b_i+b_j=b-1}}} 
   \binom a {a_i} \binom {b-1}
   {b_j} i R_i(a_i,b_i) R_j(a_j,b_j) \\
& & + 4 \binom b 2  \sum_{ \substack{{a_i + a_j=a+2}\\{b_i+b_j=b-2}}} 
  \binom a {a_i-1}
   \binom {b-2} {b_j} R_i(a_i,b_i) R_j(a_j,b_j)
  \Biggr].
\end{eqnarray*}
In each sum, it is assumed that $i,j>0$; $a_i$, $a_j$, $b_i$, $b_j
\geq 0$; $a_i + b_i = 3i-1$; $a_j + b_j = 3j-1$; and that all of these
are integers.  The large bracket corresponds to maps from reducible
curves.  The first sum in the large bracket corresponds to the case
where no tangent lines pass through the image of the node; the second
sum corresponds to when one tangent line passes through the image of
the node; and the third to when two tangent lines pass through the
image of the node (see Lemma \ref{ginger}).  Note that in the second
sum, $3i-1$ of the $a+b$ conditions fix the component corresponding to
$R_i$ (up to a finite number of possibilities).  The component
corresponding to $R_j$ is specified by the remaining $3j-2$
conditions, plus the condition that it intersect the other component
on a fixed line.

This completes the computation of the characteristic numbers for
rational plane curves.

Pandharipande earlier obtained (by topological recursion methods and 
descendants) what can be seen to be the same
recursion in the form of a differential equation (\cite{pl}):  if
$$
R(x,y,z) = \sum_{a,b,d} R_d(a,b) \frac{x^a } {a!} \frac {y^b} {b!} e^{dz},
$$
then
$$
R_{yz} = 
- R_x + R_{xz} - \frac 1 2 R^2_{zz} + (R_{zz} + y R_{xz})^2.
$$

A similar argument applied to the one-parameter family corresponding
to degree $d$ rational curves with a node at a fixed point, through
$a$ general points and tangent to $b$ general lines (where $a+b=3d-4$)
gives the formula shown in Appendix \ref{coconut}.  The corresponding
differential equation is:
\begin{eqnarray*}
NP_{yz} &=& - NP_x + NP_{xz} - \frac 1 2 R_{zzx}^2 + (R_{zzx}+ yR_{zxx})^2 \\
& & + 2 (R_{zz} + y R_{zx}) (NP_{zz} + y NP_{zx}) - R_{zz} NP_{zz}.
\end{eqnarray*}

\subsection{Swapping incidences for tangencies:  genus 1}
\label{minnow}
On the universal curve over $\mos$, let $Q$ be the divisor
corresponding to nodal irreducible fibers.  Following \cite{bbcII}, let $R$ be
the divisor corresponding to rational components of reducible fibers.  Then 
\begin{equation}
\label{mrshowell}
\om \cong \frac Q {12} + R
\end{equation}
(Kodaira's formula for the canonical bundle of an elliptic surface; see 
\cite{ccs} Theorem 12.1 for a proof over $\com$).  Hence $B =
\pi_* (D
\cdot \om) = \frac d {12} \De + \sum_i i \De_{0,i}$, so
\begin{equation}
\label{TL1}
TL = A +  \frac d {12} \De + \sum_i i \De_{0,i}.
\end{equation}

Restricting this identity to the one-parameter family corresponding to
degree $d$ elliptic curves through $a$ general points and tangent to
$b$ general lines (where $a+b=3d-1$) gives:
\begin{eqnarray*}
E_d(a,b+1) &=& E_d(a+1,b)  \\
& & + \frac d {12} \left( \binom {d-1} 2 R_d(a,b) +
2b NL_d(a,b-1) + 4 \binom b 2 NP_d(a,b-2) \right) \\
& & + \sum_{i+j=d} i \Biggl[ 
    \sum_{\substack{{a_i + a_j=a}\\{b_i+b_j=b}}} \binom a {a_i} \binom b 
   {b_i}( ij) R_i(a_i,b_i) E_j(a_j,b_j) \\
& & + 2b \biggl( \sum_{\substack{{a_i + a_j=a+1}\\{b_i+b_j=b-1}}} 
   \binom a {a_j} 
  \binom   {b-1} {b_i} j R_i(a_i,b_i) E_j(a_j,b_j) \\
& & +  \sum_{\substack{{a_i + a_j=a+1}\\{b_i+b_j=b-1}}} \binom a {a_i} \binom 
  {b-1} {b_i} i R_i(a_i,b_i) E_j(a_j,b_j) \biggr) \\
& & + 4 \binom b 2
  \sum_{\substack{{a_i + a_j=a+2}\\{b_i+b_j=b-2}}} \binom a {a_i-1} \binom 
  {b-2} {b_i}  R_i(a_i,b_i) E_j(a_j,b_j) \Biggr].
\end{eqnarray*}

$NL_d(a,b-1)$ can be found using (\ref{nl}).  The large square bracket
corresponds to maps of reducible curves.  The first sum corresponds to
the case when no tangent line passes through the image of the node,
the next two sums correspond to when one tangent line passes through
the image of the node, and the last sum corresponds to when two
tangent lines pass through the image of the node.

The corresponding differential equation is:
$$
E_y = E_x + \De + 2(R_{zz} + R_{zx}) (E_z + E_x) - R_{zz}E_z
$$
where
$$
\De = \frac 1 {12} \left( \frac 1 2 ( R_{zzz} - 3 R_{zz} + 2 R_z )
+ 2 y NL_z + 2 y^2 NP_z  \right).
$$

This completes the computation of the characteristic numbers of
elliptic plane curves.

\subsection{Characteristic numbers of elliptic curves with fixed 
$j$-invariant ($j \neq \infty$)} \label{jsub} Let $M_j$ be the Weil
divisor on $\mos$ corresponding to curves whose stable model has fixed
$j$-invariant $j$.  Then $M_j \cong M_{\infty}$ if $j \neq 0, 1728$,
$M_0 \cong M_{\infty} / 3$, and $M_{1728} \cong M_{\infty}/2$
(\cite{pj} Lemma 4).  If $a+b=3d-1$, define $J_d(a,b) := M_{\infty}
A^a TL^b$.  By Corollary \ref{enum}, if $d \geq 3$, the characteristic
numbers of curves with fixed $j$-invariant $j \neq 0, 1728, \infty$
are given by $J_d(a,b)$, and if $j=0$ or $j=1728$, the the
characteristic numbers are one third and one half $J_d(a,b)$
respectively.  But $M_{\infty}$ parametrizes maps from nodal rational
curves, so we can calculate $M_{\infty} A^a TL^b$ using Lemma
\ref{ginger}: $$ J_d(a,b) = \binom {d-1} 2 R_d(a,b) + 2b NL_d(a,b-1) +
4 \binom b 2 NP_d(a,b-2).$$

\subsection{Numbers}

Using the recursions given above, we find the following
characteristic numbers for elliptic curves.  (The first number in each
sequence is the number with only incidence conditions; the last is the
number with only tangency conditions.)

Conics:  0, 0, 0, 0, 2, 10, 45/2.

Cubics:  1, 4, 16, 64, 256, 976, 3424, 9766, 21004, 33616.

Quartics: 225,  1010, 4396, 18432, 73920, 280560, 994320, 3230956, 9409052,
23771160, 50569520, 89120080, 129996216.

Quintics: 87192, 411376, 1873388, 8197344, 34294992, 136396752, 512271756,
1802742368, 5889847264, 17668868832, 48034104112, 116575540736,
248984451648, 463227482784, 747546215472, 1048687299072.

The cubic numbers agree with those found by Aluffi in \cite{asmooth}.
The quartic numbers agree with the predictions of Zeuthen (see
\cite{s} p. 187).

Using the recursion of Subsection \ref{jsub}, we find the following
characteristic numbers for elliptic curves with fixed $j$-invariant
($j \neq 0, 1728, \infty$).

Conics:  0, 0, 0, 12, 48, 75.

Cubics:  12, 48, 192, 768, 2784, 8832, 21828, 39072, 50448.

Quartics: 1860, 8088, 33792, 134208, 497952, 1696320, 5193768,
13954512, 31849968, 60019872, 92165280, 115892448.

The cubic numbers agree with those found by Aluffi in \cite{aj} Theorem
III(2).  The incidence-only numbers necessarily agree with the numbers
found by Pandharipande in \cite{pj}, as the formula is the same.

\subsection{Characteristic numbers in $\proj^n$}
The same method gives a program to recursively compute characteristic
numbers of elliptic curves in $\proj^n$ that should be simpler than
the algorithm of \cite{char}.  Use Kontsevich's cross-ratio method to
count irreducible nodal rational curves through various linear spaces
and where the node is required to lie on a given linear space
(analogous to the derivation of (\ref{mrhowell})).  Use (\ref{TL0}) to
compute all the characteristic numbers of each of these families of
rational curves.  Use \cite{ratell} to compute the number of elliptic
curves through various linear spaces.  Finally, use (\ref{TL1}) to
compute all characteristic numbers of curves in $\proj^n$. The same
calculations also allow one to compute characteristic numbers of
elliptic curves in $\proj^n$ with fixed $j$-invariant.

\subsection{Covers of $\proj^1$}  
\label{covers}
By restricting Pandharipande's relation (\ref{TL0}) and relation
(\ref{TL1}) to degree $d$ covers of a line by a genus 0 and 1 curve
respectively (so $A$ restricts to 0), where all but 1 ramification are
fixed, we obtain recursions for $M^g_d$ ($g=0,1$), the
number of distinct covers of $\proj^1$ by irreducible genus $g$ curves
with $2d+2g-2$ fixed ramification points:
\begin{eqnarray*}
M^0_d &=& \frac {(2d-3)} d \sum_{j=1}^{d-1} \binom {2d-4} {2j-2} M^0_j M^0_{d-j} j^2 (d-j)^2 \\
M^1_d &=& \frac d 6 \binom d 2 (2d-1) M^0_d + \sum_{j=1}^{d-2} 2j (2d-1) 
\binom {2d-2}{2j-2} M^0_j M^1_{d-j} (d-j) j.
\end{eqnarray*}
The first equation was found earlier by Pandharipande and the second
by Pandharipande and Graber (\cite{gp2}); their proofs used an:
analysis of the divisors on $\mbar_{g,n}(\proj^1,d)$.  The closed-form
expression $M^0_d = d^{d-3} (2d-2)! / d!$ follows by an easy
combinatorial argument from the first equation using Cayley's formula
for the number of trees on $n$ vertices.  (This formula was first
proved in \cite{ct}.  A more general formula was stated by Hurwitz and
was first proved in \cite{gj}.  For more on this problem, including
history, see \cite{gl}.)

Graber and Pandharipande have conjectured a similar formula for $g=2$:
\begin{eqnarray*}
M^2_d  &=&  d^2 \left( \frac {97}{136}  d - \frac {20}{17} \right) M^1_d
   + \sum_{j=1}^{d-1} M^0_j M^2_{d-j} \binom {2d} {2j-2} j(d-j)
                                   \left( -\frac{115}{17} j + 8d \right)\\
  & &  + \sum_{j=1}^{d-1} M^1_j M^1_{d-j} \binom {2d}{2j} j(d-j)
              \left( \frac{11697}{34} j(d-j) - \frac{3899}{68} d^2 \right).
\end{eqnarray*}
It is still unclear why a genus 2 relation should exist (either
combinatorially or algebro-geometrically).  The relation looks as though
it is induced by a relation in the Picard group of the moduli space, but
no such relation exists.

\subsection{Divisor theory on $\cmbar_1(\proj^2,d)^*$}
In \cite{iqd}, Pandharipande determined the divisor theory on
$\cmbar_0(\proj^n,d)$ (including the top intersection products of
divisors).  The divisor theory of $\cmbar_1(\proj^2,d)^*$ is more
complicated.  In addition to the divisor $A$ and the enumeratively
meaningful boundary divisors, there are potentially three other
enumeratively meaningless divisors (see \cite{ratell} Lemma 3.14):
\begin{enumerate}
\item points corresponding to cuspidal rational curves with a contracted 
elliptic tail,
\item points corresponding to a contracted elliptic component attached to 
two rational components, where the images of the rational components meet 
at a tacnode, and
\item points corresponding to contracted elliptic components attached to 
three rational components.
\end{enumerate}
The stack $\mos$ is smooth away from these divisors.
$\cmbar_1(\proj^2,d)$ is unibranch at the third type of divisor;
Thaddeus has shown that $\mos$ is singular there (\cite{thads}).
There are several natural questions to ask about the geometry and
topology of $\mos$.  Is it smooth at the other two divisors?  Is the
normalization of $\mos$ smooth?  If $d=3$, how does it compare to
Aluffi's space of complete cubics?  What are the top intersection
products of these divisors?  (The arguments here allow us to calculate
$A^a B^{3d-a}$ and $A^a B^{3d-1-a} D$ where $D$ is any boundary
divisor.)  What about $\cmbar_1(\proj^n,d)^*$?

\section{``Codimension 1'' Numbers}

Fix a degree $d$ and geometric genus $g$.  In \cite{dh1}, Diaz and Harris
express over twenty divisors on the normalization of the Severi
variety as linear combinations of $A$, $B$, $C$, and boundaries
$\De_0$ and $\De_{i,j}$ (and conjecture that all divisors are
linear combinations).  For example, if $CU$ is the divisor of cuspidal
curves, then $CU=3A+3B+C-\De$ (\cite{dh1} (1.1)).  If $K_W$ is the
canonical bundle of the (normalization of the) Severi variety, then
$K_W = -3 A/2 + 3B/2 + 11C/12 - 13 \De/12$ (\cite{dh1} (1.17)).

Restricting these divisors to the one-dimensional family of geometric
genus $g$ degree $d$ plane curves through $3d+g-2$ general points
(which misses the enumeratively meaningless divisors), we obtain
recursive equations for the number of such curves with various
geometric behaviors (e.g. with a tacnode, three collinear nodes,
etc.).  We will give examples from the literature that turn out to be
immediate consequences of \cite{dh1}.

\subsection{Geometric and arithmetic sectional genera of the Severi variety}

We also obtain recursions for versions of the geometric and arithmetic
sectional genera.  Following \cite{canonical} Section 3, consider the curves $C_d$
(the intersection of the Severi variety with $3d+g-2$ hyperplanes
corresponding to requiring the curve to pass through $3d+g-2$ general
points $p_1$,
\dots, $p_{3d+g-2}$), $\hC_d$ (the one-parameter family of
$\cmbar_g(\proj^2,d)$ corresponding to requiring the image curve to
pass through $3d+g-2$ general points), and $\tC_d$ (the normalization
of $\hC_d$).  Let the arithmetic genera of these curves be $g_d$,
$\hg_d$, and $\tg_d$ respectively.  There are natural maps $\tC_d
\rightarrow \hC_d \rightarrow C_d$.  The singularities of $\hC_d$ are
simple nodes, which occur when the image curve has a simple node at
one of the general points $p_i$ (\cite{canonical} Section 3; the
argument holds for any $g$).  The singularities of $C_d$ are the
above, plus simple cusps corresponding to cuspidal curves, plus
singularities of the type of the coordinate axes at the origin in
$\com^{ij}$ corresponding to curves with two components (of degrees
$i$, $j$) whose geometric genera add to $g$, plus the singularities of
the type of the coordinate axes in $\com^{\binom{d-1} 2 - (g-1)}$,
corresponding to irreducible curves of geometric genus $g-1$
(\cite{dh1} Section 1).

Thus
\begin{eqnarray}
g_d - \hg_d = CU_{d,g} & +& \frac 1 2 \sum_{\substack{{i+j=d} \\ {g_i+g_j=g}}} 
(ij-1) \binom {3d+g-2} {3i+g_i-1} N^{i,g_i} N^{j,g_j} \nonumber \\
&+& \left( \binom {d-1} 2 - g \right) N^{d,g-1},
\label{gilligan}
\end{eqnarray}
where $CU_{d,g}$ is the number of irreducible degree $d$ geometric
genus $g$ cuspidal curves through $3d+g-2$ fixed general points, and
$N^{d,g}$ is the number of irreducible degree $d$ geometric genus $g$
curves through $3d+g-1$ points.

Also, Pandharipande's genus 0 argument of \cite{canonical}
3.4 works for any genus, and shows that 
\begin{equation}
\label{skipper}
\hg_d - \tg_d = (3d+g-2) NP_{d,g}
\end{equation}
where $NP_{d,g}$ is the number of irreducible degree $d$
geometric genus $g$ plane curves through $3d+g-3$ fixed general points
with a node at another fixed point.

The arithmetic (resp. geometric) sectional genus of a variety $V
\subset \proj^n$ of dimension $e$ is defined to be the arithmetic
(resp. geometric) genus of the curve obtained by intersecting $V$ with
$e-1$ general hyperplanes.

\begin{pr}
The geometric sectional genus is $\hg_d$.
\end{pr}
\bpf
From \cite{dh1} Section 1, the only codimension 1 singularities of the
Severi variety $V^{d,g}$ are those corresponding (generically) to
cuspidal curves and curves with $\de+1$ nodes ($\de := \binom{d-1} 2 -
g$), and the singularities are as described above.  If $V^{d,g}$ is
intersected with (special) hyperplanes corresponding to requiring the
curve to pass through various generally chosen fixed points, the
intersection picks up new singularities, corresponding to curves with
a node at one of the fixed points.  Hence the geometric sectional
genus is the genus of the partial normalization of $C_d$ corresponding
to normalizing the singularities corresponding to cuspidal and
$(\de+1)$-nodal curves, which is the arithmetic genus of $\hC_d$.
\epf

{\em Notational caution:} In \cite{canonical}, $\hg_d$ is called the ``arithmetic
genus''.

\subsection{Genus 0}

Clearly, $|A| = R_d$.  By (\ref{TL0}), 
$$
|B| = - \frac {R_d} d + \frac 1 {2d} \sum_{i+j=d} \binom {3d-2}{3i-1} i^2 j^2 R_i R_j.
$$
It is simple to show (e.g. \cite{iqd} Lemma 2.1.2) that $C=-\De$, so
$$
|C| = - \frac 1 2 \sum_{i+j=d} \binom{3d-2}{3i-1} ij R_i R_j.
$$
Note that Kontsevich's recursion (\ref{Kont}) can be rewritten as
\begin{equation}
\label{Kont2}
9(d-2)A = 3(d+2)B + 2dC
\end{equation}
(or $\pi_* (3D + \om) \cdot (3 (d-2) D - 2d \om)$ restricted to the one-parameter family is 0).

The formula of Katz-Qin-Ruan for the number of degree $d$
triple-pointed rational curves (\cite{kqr}, Lemma 3.2) can be
rewritten as
\begin{equation} \label{triple}
(d^2-6d+10)|A|/2 - (d-6)|B|/2+|C|
\end{equation}
which is the $g=0$ case of \cite{dh1} (1.3).  Pandharipande's formula
for the number of degree $d$ rational cuspidal curves (\cite{iqd}
Prop. 5) can be rewritten as $3|A| + 3 |B| + |C| - |\De|$, which is
the $g=0$ case of \cite{dh1} (1.1).  Ran's formula for the cuspidal
number (\cite{bbcII} Theorem (ii) (2)) yields the same numbers for
small $d$, and presumably is the same formula after a substitution.

By adjunction, the geometric sectional genus $\hg_d$ of the Severi
variety is given by $2 \hg_d - 2 = | K_W + (3d-2)A |$.  The formula of
Pandharipande for $\hg_d$ (\cite{canonical} Section 3.2) can be rewritten as $$ 2
\hat{g}_d - 2 = ( -3 |A|/2 + 3 |B|/2 + 11 |C| / 12 - 13 |\De|/12) +
(3d-2) |A|, $$
 which is the $g=0$ case of \cite{dh1} (1.17).
(Pandharipande then computes the arithmetic sectional genus $g_d$ 
using (\ref{gilligan}).  His
computation of $\tg_d$ by other means gives his recursive formula for
$NP_d$ (mentioned in Subsection \ref{theprofessor}) via (\ref{skipper}).)

\subsection{Genus 1}
Clearly $|A|= E_d$ and
$$
|\De| = \binom {d-1} 2 R_d + \sum_{i+j=d} ij \binom {3d-1} {3i-1} R_i E_j.
$$
From Subsection \ref{minnow}, $B=\frac d {12} \De_0 + \sum_i i \De_{0,i}$, so
$$
|B| = \frac d {12} \binom {d-1} 2 R_d + \sum_{i+j=d} i^2 j \binom {3d-1} {3i-1} R_i E_j.
$$
From the description of $\om$ in Subsection \ref{minnow}, 
$$
|C| = -  \sum_{i+j=d} ij \binom {3d-1} {3i-1} R_i E_j.
$$
Note that the Eguchi-Hori-Xiong recursion can be rewritten as $9A-3B-2C=0$ 
(or $\pi_* (3D + \om) \cdot (3 D - 2 \om)$ restricted to the one-parameter 
family is numerically 0, cf. (\ref{Kont2})).

Ran's formula for the number of degree $d$ cuspidal elliptic curves
(\cite{bbcII} Theorem (ii) (3)) can be rewritten as $|3A+3B+C - \De|$, which is
the $g=1$ case of \cite{dh1} (1.1).  Call this number $CU_{d,1}$.

Using \cite{dh1} as in the genus 0 case, we find the geometric sectional
genus of the Severi variety $\hg_d$:
\begin{eqnarray*}
 2 \hat{g}_d - 2 &=& ( -3 |A|/2 + 3 |B|/2 + 11 |C| / 12 - 13 |\De|/12) +
(3d-1) |A|, \\
&=& \left( 3d - \frac 5 2 \right) E_d + \left( \frac {3d-26} {24} \right) 
   \binom{d-1} 2 R_d  \\
& & + \sum_{i+j=d} ij \binom{3d-1}{3i-1} R_iE_j \left( \frac 3 2 i-2 \right).
\end{eqnarray*}
This formula is identical to that of Ran's Theorem (ii) of \cite{bbcII}.  Via
(\ref{gilligan}), this yields a recursion for the arithmetic
sectional genus of the Severi variety$g_d$: $$ g_d = \hg_d + CU_{d,1}
+ \sum_{i+j=d} (ij-1) \binom{3d-1}{3i-1} R_i E_j + \left( \binom {d-1}
2 - 1 \right) R_d.  $$

The values of $\hg_d$ for $3 \leq d \leq 7$ are: 0, 486, 410439,
395296561, 534578574561.  The values of $g_d$ for $3 \leq d \leq 7$
are: 0, 2676, 1440874, 1117718773, 1317320595961.

\subsection{Genus 2}

Let $T_d$ be the number of irreducible degree $d$ geometric genus 2
plane curves through $3d+1$ fixed general points ($d>2$).  From
\cite{rinv} or \cite{ch}, the numbers $|A|$ and $|B|$ can be found
(the latter by computing $|TL| = |A| + |B|$, the number of irreducible
geometric genus 2 plane curves through $3d$ points tangent to a fixed
line).  The number $|A|$ can be computed more easily by the recursion
of Belorousski and Pandharipande \cite{bp}.  (Their ideas should also
lead to a recursive calculation for $|B|$.)  Also, $$ |\De| = \left(
\binom{d-1} 2 - 1 \right) E_d +
\sum_{i+j=d} ij \left( \binom {3d} {3i-1} R_i T_j + \frac 1 2 \binom
{3d} {3i} E_i E_j \right).  $$

To compute $|C|$, consider the family of genus 2 curves to be pulled
back from the universal curve over the moduli stack $\cmbar_2$, blown
up at a finite number of points (corresponding to the points in the
family where the curve is a genus 2 curve and a genus 0 curve
intersecting at a node).  If $\rho: \cu \rightarrow \cmbar_2$ is the
universal curve over $\cmbar_2$, and $\om_{\rho}$ is the relative
dualizing sheaf, then by \cite{m} (8.5), $$
\rho_* ( \om_{\rho}^2) = (\de_0 + 7 \de_1)/5
$$
where $\de_0$ is the divisor corresponding irreducible nodal curves
and $\de_1$ is the divisor corresponding to reducible nodal curves
(with each component of genus 1).  Hence $|C|$ can be expressed in terms 
of previously-known quantities:
$$
|C| = \frac 1 5 \left( \binom {d-1} 2 - 1 \right) E_d + \frac 7 {10}
\sum_{i+j=d} ij \binom {3d} {3i} E_i E_j - \sum_{i+j=d} i j \binom {3d} {3i-1} R_i T_j.
$$
Examples are given at the end of the section.

\subsection{Genus 3}  
Once again, $|A|$ and $|B|$ can be calculated by the algorithm of
\cite{rinv} or \cite{ch}, and $|\De|$ can be inductively calculated.
Graber has found a recursive method of counting the number of genus
$g$ hyperelliptic plane curves through $3d+1$ general points
(\cite{g}) by relating these numbers to the Gromov-Witten invariants
of the Hilbert scheme of two points in the plane.  (The algorithm is
effective, and maple code is available.)  Call the genus 3
hyperelliptic numbers $H_d$; the smallest non-zero values are $H_5 =
135$, $H_6=3929499$, $H_7=23875461099$ (\cite{g}).  If $h$ is the
reduced divisor of the hyperelliptic locus on the stack $\cmbar_3$,
then $h=9 \la - \de_0 - 3
\de_1$ (see \cite{h} appendix for explanation and proof).  As in the
genus 2 case, if $\rho$ is the structure map of the universal curve
over $\cmbar_3$, $\rho_* (\om_{\rho}^2) = 12 \la - \de_0 - \de_1$ (see
\cite{m} p. 306), so $\rho_* (\om_{\rho}^2) = (4 h + \de_0 + 9 \de_1)/3$.  
Hence
\begin{eqnarray*}
|C| &= &\frac 4 3 H_d + \frac 1 3 \left( \binom {d-1} 2 - 2 \right) T_d \\
& & +  \sum_{i+j=d} ij \left( 3 \binom {3d+1} {3i} E_i T_j -
\binom {3d+1} {3i-1} R_i U_j \right)
\end{eqnarray*}

In this way, all codimension 1 numbers for genus 2 and 3 curves can be
computed.  As examples, for $4 \leq d \leq 6$, $|A|$, $|B|$, $|C|$,
$|\De|$, and $|TL|$ are given as well as $|CU|$, the number of
cuspidal curves, and $\hg$ and $g$, the geometric and arithmetic
sectional genera of the Severi variety.

\begin{tabular}{|c|c|c|c|c|c|c|} \cline{2-7}
\multicolumn{1}{c}{} 
& \multicolumn{3}{|c|} {$g=2$} & \multicolumn{3}{|c|} {$g=3$}  \\ 
\cline{2-7} \multicolumn{1}{c}{} & \multicolumn{1}{|c|}{$d=4$} 
        & \multicolumn{1}{c|}{$d=5$} 
        & \multicolumn{1}{c|}{$d=6$} 
        & \multicolumn{1}{c|}{$d=4$} 
        & \multicolumn{1}{c|}{$d=5$} 
        & \multicolumn{1}{c|}{$d=6$} \\ \hline
$|A|$   & 27    & 36855    &  58444767 & 1 & 7915 & 34435125 \\
$|B|$   & 117   & 166761   & 268149471 & 5 & 41665 & 182133909 \\
$|C|$   & 90    &  75852   &  73644975 & 9 & 48840 & 154231695 \\
$|\De|$ & 450   & 447300   & 547180713 & 27 & 147900 & 474418485 \\
$|TL|$  & 144   & 203616   & 326594238 & 6 & 49580 & 216569034 \\
$|CU|$  & 72    & 239400   & 506246976 & 0 & 49680 & 329520312 \\ 
$\hg$   & 28    & 166321   & 420645826 & 0 & 30906 & 251620624 \\ 
$g$     & 325   & 762994   & 1410743814 & 0 & 191511 & 995749561 \\ \hline
\end{tabular}

\appendix
\section{A recursive formula for $NP(a,b)$}
\label{coconut}
\begin{eqnarray*}
NP(a,b+1) &=& \frac {d-1} d NP(a+1,b)  \\
& & + \sum_{i+j=d} \frac {ij} {2d} \Biggl[
\sum_{\substack{{a_i + a_j=a+2}\\{b_i+b_j=b}}} \binom a {a_i-1} \binom b {b_i} (ij-1) 
   R_i(a_i,b_i) R_j(a_j,b_j) \\
& & + 2 \sum_{\substack{{a_i + a_j=a}\\{b_i+b_j=b}}} \binom a {a_i} \binom b
   {b_i} (ij) R_i(a_i,b_i) NP_j(a_j,b_j) \\
& & + 4b \sum_{\substack{{a_i + a_j=a+3}\\{b_i+b_j=b-1}}} \binom a {a_i-1} \binom {b-1}
   {b_i} i R_i(a_i,b_i) R_j(a_j,b_j) \\
& & + 4b \sum_{\substack{{a_i + a_j=a+1}\\{b_i+b_j=b-1}}} \binom a {a_i} \binom {b-1}
   {b_i} i NP_i(a_i,b_i) R_j(a_j,b_j) \\
& & + 4b \sum_{\substack{{a_i + a_j=a+1}\\{b_i+b_j=b-1}}} \binom a {a_i} \binom {b-1}
   {b_i} i R_i(a_i,b_i) NP_j(a_j,b_j) \\
& & + 4 \binom b 2 \sum_{\substack{{a_i + a_j=a+4}\\{b_i+b_j=b-2}}} \binom a {a_i-2}
   \binom {b-2}  {b_i} R_i(a_i,b_i) R_j(a_j,b_j) \\
& & + 8 \binom b 2 \sum_{\substack{{a_i + a_j=a+2}\\{b_i+b_j=b-2}}} \binom a {a_i-1}
   \binom {b-2}  {b_i} R_i(a_i,b_i) NP_j(a_j,b_j)  \Biggr].
\end{eqnarray*}

In each sum in the large bracket, it is assumed that $a_i+b_i=3i-1$ if
$R_i(a_i,b_i)$ appears in the sum, and $a_i+b_i=3i-3$ if $NP_i(a_i,b_i)$
appears.  The same assumption is made when $i$ is replaced by $j$.

The large square bracket corresponds to maps from reducible curves.
(To avoid confusion: the ``image of the node'' refers to the image of
the node of the source curve.  The ``fixed node'' refers to the node
of the {\em image} that is required to be at a fixed point.)  Zero,
one, or two tangent lines can pass through the image of the node of
the source curve.  The two branches through the fixed node can belong
to the same component, or one can belong to each.  The table below
identifies which possibilities correspond to which sum in the large bracket.

\begin{tabular}{|l|c|c|} \hline
sum  & number of tangent   & number of \\
 & lines through image & irreducible components \\
 & of node of source &  through fixed node \\
 \hline
first  & 0  & 2 \\
second  &  0 & 1\\
third   &1 &2\\
fourth and & & \\
 \; fifth  &  1 & 1\\
sixth  &2 &2 \\
seventh  &2 &1 \\ \hline
\end{tabular}

\end{document}